\documentstyle[11pt,amsfonts,amscd]{amsproc}
\newtheorem{Theorem}{Theorem}[section]
\newtheorem{Definition}[Theorem]{Definition}
\newtheorem{Proposition}[Theorem]{Proposition}

\newtheorem{Lemma}[Theorem]{Lemma}
\newtheorem{Corollary}[Theorem]{Corollary}
\theoremstyle{remark}
\newtheorem{Example}[Theorem]{Example}

\def\eps{\varepsilon}

\def\ovr{\overline}

\def\om{\omega}

\def\al{\alpha}
\def\gm{\gamma}
\def\th{\theta}

\def\dl{\delta}

\def\lm{\lambda}

\def\si{\sigma}

\def\sbs{\subset}

\def\be{\begin{enumerate}}
\def\ee{\end{enumerate}}
\def\bT{\begin{Theorem}}
\def\eT{\end{Theorem}}
\def\bP{\begin{Proposition}}
\def\eP{\end{Proposition}}
\def\bD{\begin{Definition}}
\def\eD{\end{Definition}}
\def\bE{\begin{Example}}
\def\eE{\end{Example}}
\def\bL{\begin{Lemma}}
\def\eL{\end{Lemma}}
\def\bC{\begin{Corollary}}
\def\eC{\end{Corollary}}
\def\A{{\mathcal A}}

\def\A{{\mathcal A}}
\def\B{{\mathcal B}}

\def\M{{\mathcal M}}

\def\P{{\mathcal P}}
\def\R{{\mathcal R}}
\def\Q{{\mathcal Q}}

\begin{document}
\title{Stable algebras of entire functions}
\author{ Dan Coman and Evgeny A. Poletsky}
\thanks{Both authors are supported by NSF Grants.}
\subjclass[2000]{ Primary: 32A38. Secondary: 30H05}
\address{ Department of Mathematics,  215 Carnegie Hall,
Syracuse University,  Syracuse, NY 13244-1150, USA. E-mail:
dcoman@@syr.edu, eapolets@@syr.edu}
\begin{abstract} Suppose that $h$ and $g$ belong to the algebra
$\B$ generated by the rational functions and an entire function
$f$ of finite order on ${\Bbb C}^n$ and that $h/g$ has algebraic
polar variety. We show that either $h/g\in\B$ or $f=q_1e^p+q_2$,
where $p$ is a polynomial and $q_1,q_2$ are rational functions. In
the latter case, $h/g$ belongs to the algebra generated  by the
rational functions, $e^p$ and $e^{-p}$.
\end{abstract}
\maketitle
\section{Introduction}
\par In this paper we approach the problem of algebraic
independence of entire functions on ${\Bbb C}^n$ over ${\Bbb C}$
or over the ring of polynomials $\P^n$ on ${\Bbb C}^n$. Two entire
functions $f$ and $g$ are algebraically dependent if there is a
non-zero polynomial $P$ in ${\Bbb C}[x,y]$ or in $\P^n[x,y]$ such
that $P(f,g)\equiv 0$.
\par The first case, i.e., $P\in{\Bbb C}[x,y]$, can be answered
completely and this is done in Section \ref{S:pairs}. We prove in
Theorem \ref{T:entpair} that if $f$ and $g$ are algebraically
dependent entire functions over ${\Bbb C}$, then there is an
entire function $h$ such that either $f,g\in{\Bbb C}[h]$ or
$f,g\in{\Bbb C}[e^h,e^{-h}]$. In Theorem \ref{T:meropair} such a
description is given for meromorphic functions $f$ and $g$.
\par To approach the seemingly difficult problem of algebraic
dependence of functions over $\P^n$ we study in this paper the
dependence relation $P(f,g)=P_1(f)g+P_0(f)\equiv0$, where
$P_1,P_0\in\P^n[x]$. In other words, if $f$ is entire we are
asking when the ratio $P_0(f)/P_1(f)$ is also entire.
\par This problem has a rather long history. Following \cite{BD},
we say that a subalgebra $\B$ of an algebra $\A$ is {\it stable
in} $\A$ if whenever $g,h\in\B$ and $h/g\in\A$ then $h/g\in\B$.
\par In 1929 Ritt (\cite{Ri1},\cite{Ri2}) proved that the
algebra of quasipolynomials $\Q$ is stable in the algebra of
entire functions on ${\Bbb C}$. Quasipolynomials are linear
combinations of exponentials $e^{\lm_jz}$. In \cite{Sh} Shields
improved Ritt's theorem by showing that $\Q$ is also stable in the
algebra of meromorphic functions whose number of poles in $|z|<r$
is $o(r)$ as $r\to\infty$. Gordon and Levin considered in
\cite{GoLe} algebras of more general quasipolynomials on angles in
${\Bbb C}$ and proved the stability property in these cases.
\par In \cite{BD} Berenstein and Dostal considered the field $\R^n$
of rational functions on  ${\Bbb C}^n$ and the algebra of entire
$\R^n$-quasipolynomials. These are functions of the form
$\sum_{j=1}^kr_j(z)e^{\langle\theta_j,z\rangle}$, where
$\theta_j\in{\Bbb C}^n$ and $r_j\in\R^n$. They proved that this
algebra is stable in the algebra of all entire functions. The
choice of $\R^n$ over $\P^n$ was done to accommodate the example
of the function $z^{-1}\sin z$ on ${\Bbb C}$.
\par One may think that the stability of algebras of
quasipolynomials is due to the special features of exponentials.
However, in this paper we show that in fact the stability is a
generic property of algebras of entire functions. Let $\M_0^n$ be
the algebra of meromorphic functions on ${\Bbb C}^n$ whose polar
varieties are algebraic. We prove the following theorem:
\bT\label{T:it} Let $\R^n[f]\subset\M_0^n$ be the subalgebra
generated by $\R^n$ and an entire function $f$ of finite order on
${\Bbb C}^n$. Then either $\R^n[f]$ is stable in $\M_0^n$ or
$\R^n[f]=\R^n[e^p]$ for some $p\in\P^n$ and the algebra
$\R^n[e^p,e^{-p}]$ is stable in $\M_0^n$.
\eT
\par This theorem is proved in Section \ref{S:dt}. For this
we show in Theorem \ref{T:fn0} and Corollary \ref{C:eq} that an
algebra $\R^n[f]$ contains non-trivial (i.e., non-rational)
invertible elements in $\M_0^n$ if and only if $f=q_1e^p+q_2$,
where $q_1,q_2\in\R^n$ and $p\in\P^n$.
\par The proof of our basic tool, namely Theorem \ref{T:fn0}, relies
on solving  in Section \ref{S:aeina} algebraic equations in an
angle in ${\Bbb C}$. Let $P$ be a monic polynomial with slowly
growing holomorphic coefficients in an angle. We show in Theorem
\ref{T:mt} that under some suitable growth conditions on $e^p$ the
equation $P(f)=e^p$ has a holomorphic solution if and only if
$P(w)=(w+q)^m$.
\par The results in Section \ref{S:aeina} use Phragmen--Lindel\"of type
theorems presented in Section \ref{S:pr}. They relate the growth
of holomorphic functions in an angle with the behavior of their
indicators.
\section{Preliminary results}\label{S:pr}
\par We will need the following standard estimate for the roots of a
polynomial. \bP\label{P:er} If
$P(z)=z^m+a_{m-1}z^{m-1}+\dots+a_1z+a_0=0$ and
$\|P\|=\max\{1,|a_0|,\dots,|a_{m-1}|\}$ then
$|a_0|(m\|P\|)^{-1}\le |z|\le m\|P\|$.
\eP
\par For $\al<\beta$ and $r_0\ge0$ let
$$S(\al,\beta,r_0)=\{z=re^{i\th}:\,r>r_0,\;\al<\th<\beta\}.$$ We say
that a holomorphic function $f$ on $S(\al,\beta,r_0)$ is {\it of
order at most $\rho\ge0$} if $|f(z)|\le C(\rho')e^{|z|^{\rho'}}$
for every $\rho'>\rho$. For every $\th$, $\al<\th<\beta$, the
indicator of $f$ of order $\rho$ is defined by
$$h(\th)=h_{f,\rho}(\th)=
\limsup_{r\to\infty}\frac{\log|f(re^{i\th})|}{r^\rho}.$$ The
indicator may have infinite values. If $h_{f,\rho}(\th)<\infty$
for $\al<\th<\beta$, then we say that the function $f$ has {\it
finite indicator of order $\rho$}.
\par The following lemma is a version of \cite[Lemma 6, p. 52]{Le}
in our setting, and is proved in the same way: \bL\label{L:l6} Let
$f$ be a holomorphic function on $S(\al,\beta,r_0)$. Suppose that
$\al<\th_1<\th<\th_2<\beta$ and $\th_2-\th_1<\pi/\rho$. If $f$ is
of order at most $\rho>0$ and has finite indicator $h=h_{f,\rho}$,
then
$$h(\th)\le\frac{h(\th_1)\sin\rho(\th_2-\th)+
h(\th_2)\sin\rho(\th-\th_1)}{\sin\rho(\th_2-\th_1)}.$$
\eL
\par It follows from this lemma that a finite indicator is
continuous. The following lemma can be easily obtained using
elementary conformal mappings. \bL\label{L:hm} Let
$D_r=S(-\al,\al,r_0)\cap\{|z|<r\}$ and
$\gm_r=\{z=re^{i\th}:\,-\al<\th<\al\}$. If $\om_r$ is the harmonic
measure of $\gm_r$ in $D_r$ and $K$ is a compact set in
$S(-\al,\al,r_0)$, then there are positive constants $C_1$ and
$C_2$ depending on $K$ such that for all $r$ sufficiently large we
have
$$C_1r^{-\si}\le\om_r(z)\le C_2r^{-\si},\;z\in K,$$
where $\si=\pi/(2\al)$. \eL
\par Using this lemma and the Two Constants Theorem, the following
results related to the Phragmen-Lindel\"of principle can be easily
derived. \bC\label{C:ef} Let $f$ be a holomorphic function on
$\ovr S(-\al,\al,r_0)$ of order at most $\rho<\pi/(2\al)$. \par
(i) If the holomorphic function $e^{-g}$ is also of order at most
$\rho$ and $|f(z)|\le \left|e^{g(z)}\right|$ on $\partial
S(-\al,\al,r_0)$, then $|f(z)|\le \left|e^{g(z)}\right|$ in
$S(-\al,\al,r_0)$. \par (ii) If $h_{f,\rho}(\pm\al)<-\dl<0$ then
there is a constant $C>0$ such that $|f(z)|\le Ce^{-\dl|z|^\rho}$
on $\ovr S(-\al,\al,r_0)$.
\par (iii) If $h_{f,\rho}(\pm\al)\le0$ then for any $a>0$ there is a
constant $C>0$ such that $|f(z)|\le Ce^{a|z|^\rho}$ on $\ovr
S(-\al,\al,r_0)$. \eC
\par While Corollary \ref{C:ef} gives upper bounds on the absolute
values of holomorphic functions on small angles, the following
lemma shows that such functions cannot be too small on big angles.
\bL\label{L:fe0} Let $f$ be holomorphic on $S(-\al,\al,r_0)$
and let $\rho>\pi/(2\al)$. \par (i) If $|f(z)|\le Ce^{-a|z|^\rho}$
on $S(-\al,\al,r_0)$, where $a,C>0$, then $f\equiv0$. \par (ii) If
$f$ is of order at most $\rho>0$, $h_{f,\rho}(\th)\le0$ for all
$\th\in(-\al,\al)$ and $h_{f,\rho}(\th_0)<0$ for some
$\th_0\in(-\al,\al)$, then $f\equiv0$. \eL
\begin{pf} The first part follows immediately from
Lemma \ref{L:hm} and the Two Constants Theorem. To prove the
second part we note that by Lemma \ref{L:l6} $h_{f,\rho}(\th)<0$
for all $\th\in(-\al,\al)$. Choose $0<\al'<\al$ such that
$2\al'\rho>\pi$. Since $h_{f,\rho}$ is a continuous function,
there is $\dl>0$ such that $h_{f,\rho}(\th)<-\dl$ on
$[-\al',\al']$. Take $\th_1$ and $\th_2$ such that
$-\al'\le\th_1<\th_2\le\al'$ and $\th_2-\th_1<\pi/\rho$. By the
second part of Corollary \ref{C:ef} there is a constant $C>0$ such
that $|f(z)|\le Ce^{-\dl|z|^{\rho}}$ on $\ovr S(\th_1,\th_2,r_0)$.
The region $S(-\al',\al',r_0)$ can be covered by finitely many
regions $S(\th_1,\th_2,r_0)$ with $\th_2-\th_1<\pi/\rho$.
Therefore $|f(z)|\le Me^{-\dl|z|^{\rho}}$ on $S(-\al',\al',r_0)$.
By the first part of this lemma $f\equiv0$.
\end{pf}
\section{Algebraic equations in an angle}\label{S:aeina}
\par Let $\A=\A_{\al,\beta,r_0,\rho}$ be the set of functions $f$
holomorphic on $S(\al,\beta,r_0)$, of order at most $\rho$ and so
that $h_{f,\rho}(\th)\le0$ for $\al<\th<\beta$. Note that $\A$ is
an algebra. Let $\B$ be a subalgebra of $\A$. We denote by $\B[w]$
the algebra of polynomials with coefficients in $\B$.
\par Let $p$ be a holomorphic function on $S(\al,\beta,r_0)$.
We say that the indicator of order $\rho>0$ of the function $e^p$
is {\it almost sinusoidal} if $e^p$ and $e^{-p}$ are of order at
most $\rho$, $h_{e^{-p},\rho}(\th)<0$ on some interval
$(\al_1,\beta_1)\sbs(\al,\beta)$, while $h_{e^p,\rho}(\th)\le 0$
for all remaining $\th$. This definition is justified by a theorem
of Cartwright \cite{C}, which states that if $h_{e^p,\rho}$ is
{\it sinusoidal} and $\beta-\al>\pi/\rho$, then $e^p$ is of
completely regular growth and has all the above properties.
\bT\label{T:mt} Let $\rho>0$ and $\beta-\al>\pi/\rho$. Suppose
that there is a monic polynomial $P\in\B[w]$ of degree $m\ge1$ and
a function $f$ holomorphic on $S(\al,\beta,r_0)$ such that
$P(f)=e^p$. If $e^p$ has almost sinusoidal indicator of order
$\rho$, then $P(w)=(w+q)^m$ for some $q\in\B$ and $f=\eps
e^{p/m}-q$, where $\eps^m=1$. \eT
\begin{pf} The case $m=1$ is trivial. Let $m\ge2$ and
$$P(f)=P_mf^m+P_{m-1}f^{m-1}+\dots+P_0=e^{p},$$
where $P_j\in\B$, $0\le j\le m-1$, and $P_m\equiv1$. If
$f_1=f+P_{m-1}/m$ then $f_1$ verifies the equation
\begin{equation}\label{e:me1}
f_1^m+R_{m-2}f_1^{m-2}+\dots+R_0=e^p,
\end{equation}
where
$$R_k=\sum_{j=k}^m{j\choose k}(-1)^{j-k}m^{k-j}P_jP_{m-1}^{j-k}\in\B.$$
\par  The function $f_2=e^{-p/m}f_1$ satisfies the equation
\begin{equation}\label{e:epj0}
f_2^m+R_{m-2}e^{-2p/m}f_2^{m-2}+\dots+R_0e^{-p}=1.
\end{equation}
Take any $\al_1'<\beta_1'$ in $(\al_1,\beta_1)$. Then
$h_{e^{-p},\rho}(\th)<-\dl<0$ for some $\dl>0$ on
$[\al_1',\beta_1']$. By Lemma \ref{L:fe0}
$\beta_1-\al_1\le\pi/\rho$, otherwise $e^{-p}\equiv0$. Thus
$\beta'_1-\al'_1<\pi/\rho$. Corollary \ref{C:ef} shows that there
is a positive constant $C_1$ such that $\left|e^{-p(z)}\right|\le
C_1e^{-\dl|z|^\rho}$ on $\ovr S(\al_1',\beta_1',r_0)$. Since
$h_{R_j,\rho}(\th)\le0$ it follows from the same corollary that
there is a constant $C_2>0$ such that for $z\in\ovr
S(\al_1',\beta_1',r_0)$
$$\left|R_{m-2}(z)\right|\le
C_2e^{a|z|^\rho},\;\;\left|R_{m-j}(z)e^{-(j-2)p(z)/m}\right|\le
C_2e^{-a|z|^\rho},\;j\ge3,$$ where $a=\dl/(2m)$.
\par  These estimates together with Proposition \ref{P:er} and
(\ref{e:epj0}) show that the function $f_2$ is bounded on $\ovr
S(\al_1',\beta_1',r_0)$. Moreover, we see that
$$f_2^m=1+G_1e^{-2p/m},$$
where $G_1$ is holomorphic and $|G_1(z)|<C_3e^{a|z|^\rho}$ on
$\ovr S(\al_1',\beta_1',r_0)$. We conclude that there exist
$r_1\ge r_0$ and $\eps$, $\eps^m=1$, such that
\begin{equation}\label{e:epj1}
f_2=\eps+Ge^{-2p/m},
\end{equation}
where $G$ is holomorphic and
\begin{equation}\label{e:epj2}
|G(z)|<C_4e^{a|z|^\rho}
\end{equation}
on $\ovr S(\al_1',\beta_1',r_1)$. Clearly $\eps$ does not depend
on our choice of $\al'_1$ and $\beta'_1$, while the constants
$C_4$ and $\dl$ may depend on the indicated choice.
\par Multiplying (\ref{e:epj1}) by $e^{p/m}$ we obtain that
$$f_1=\eps e^{p/m}+g,\;g=e^{-p/m}G,$$
hence $g$ is analytic on $S(\al,\beta,r_0)$. Using (\ref{e:epj2})
and the fact that $h_{e^{-p},\rho}(\th)<-\dl$ on
$[\al_1',\beta_1']$, we obtain $h_{g,\rho}(\th)<0$ for
$\th\in(\al_1,\beta_1)$.
\par Letting $R_m\equiv1$, $R_{m-1}\equiv0$, and plugging
$f_1=\eps e^{p/m}+g$ into (\ref{e:me1}) we get the equation
$$g^m+Q_{m-1}g^{m-1}+\dots+Q_1g+Q_0=0,$$
where
$$Q_0=\sum_{j=0}^{m-2}\left(\eps e^{p/m}\right)^{j}R_j,\;\;
Q_k=\sum_{j=k}^m{j\choose k}\left(\eps e^{p/m}\right)^{j-k}R_j,\;
1\le k\le m-1.$$
\par By Proposition \ref{P:er} the order of $g$ does not
exceed $\rho$. If $\th\not\in(\al_1,\beta_1)$ then
$h_{e^p,\rho}(\th)\le0$. Hence $h_{Q_j,\rho}(\th)\le0$. Again
Proposition \ref{P:er} yields that $h_{g,\rho}(\th)\le0$. Thus the
function $g$ satisfies the conditions of Lemma \ref{L:fe0} on
$S(\al,\beta,r_0)$ and, therefore, $g\equiv0$ and $f_1=\eps
e^{p/m}$.
\par Plugging $f_1=\eps e^{p/m}$ into (\ref{e:me1}) and dividing by
$e^{\frac{m-2}mp}$ we obtain
$$\eps^{m-2}R_{m-2}+\eps^{m-3}R_{m-3}e^{-p/m}+
\dots+R_0e^{\frac{2-m}mp}=0.$$ This immediately implies that
$h_{R_{m-2},\rho}(\th)<0$ when $\th\in(\al_1,\beta_1)$. By Lemma
\ref{L:fe0} $R_{m-2}\equiv0$. Continuing this process we conclude
that $R_{m-2}\equiv\dots\equiv R_0\equiv0$ on $S(\al,\beta,r_0)$.
\par Thus $f=\eps e^{p/m}-P_{m-1}/m$ and $P(w)=(w+P_{m-1}/m)^m$.
\end{pf}
\bC\label{C:crg}Let $\rho>0$ and $\beta-\al>\pi/\rho$. Suppose
that there is a monic polynomial $P\in\B[w]$ of degree $m\ge1$ and
a function $f$ holomorphic on $S(\al,\beta,r_0)$ such that
$P(f)=e^p$. If the function $e^p$ is of completely regular growth
of order $\rho$, then the conclusions of Theorem \ref{T:mt} hold.
\eC
\begin{pf} Since the density of zeroes of $e^p$ is zero, by
\cite[Corollary, p. 155]{Le} $h_{e^{\pm p},\rho}(\th)=\pm
a\sin\rho(\th-\th_0)$. Hence we can choose an interval
$(\al',\beta')\sbs(\al,\beta)$, $\beta'-\al'>\pi/\rho$, where the
hypotheses of Theorem \ref{T:mt} are verified.
\end{pf}
\par Let $\A^0=\A^0_{r_0,\rho}$ be the algebra of functions
$f$ holomorphic on $\{|z|>r_0\}$, of order at most $\rho$ and with
indicator $h_{f,\rho}(\th)\le0$ for all $\th$.
\bC\label{C:ef1} Suppose that $P\in\A^0[w]$ is a monic polynomial of
degree $m\ge1$ and $f,\,p$ are holomorphic functions on
$\{|z|>r_0\}$, such that $P(f)=e^p$. If $f$ is of order at most
$\rho$ then either $p$ extends analytically at $\infty$ or the
conclusions of Theorem
\ref{T:mt} hold. \eC
\begin{pf} It follows that the function $e^p$ is of order at most
$\rho$. Thus the function $p$ has a pole at infinity of order
$n\le\rho$. If $n=0$ then $p$ extends analytically at $\infty$. If
$n\ge1$ then $e^p$ is of completely regular growth of order $n$.
So the previous corollary applies.
\end{pf}
\section{Division theorems}\label{S:dt}
\par We denote by $\P^n$ and $\R^n$ the spaces of polynomials,
respectively rational functions, on ${\Bbb C}^n$. Let $\M_0^n$ be
the set of meromorphic functions on ${\Bbb C}^n$ of the form
$h/p$, where $h$ is entire and $p\in\P^n$. Equivalently, $\M_0^n$
consists of the meromorphic functions on ${\Bbb C}^n$ whose polar
variety is algebraic. Note that $\M_0^n$ is an algebra and
$\R^n\subset\M_0^n$. Moreover, the invertible elements of $\M_0^n$
are the functions $Re^h$, where $R\in\R^n$ and $h$ is entire. If
$h$ is not a constant we will call $Re^h$ a {\it non-trivial
invertible function} in $\M_0^n$.
\par We let $\R^n[w]$ be
the algebra of polynomials in $w$ with coefficients in $\R^n$. If
$P\in\R^n[w]$, the degree of $P$ will mean the degree of $P$ in
$w$. For a function $f\in\M_0^n$ we denote by
$\R^n[f]\subset\M_0^n$ the algebra generated by $\R^n$ and $f$. We
will omit the index $n$ in our notations if $n=1$.

\bT\label{T:fn0} Let $f$ be an entire transcendental function on
${\Bbb C}^n$ of finite order. If $P(f)=Re^p$ is a non-trivial
invertible function in $\M_0^n$, where $P\in\R^n[w]$ and
$R\in\R^n$, then $p\in\P^n$, $\deg p\geq1$, and there exist an
integer $m\ge1$ and $p_1,p_2\in\R^n$ such that
$$P(w)=Rp_2^m\left(w+p_1\right)^m\;,\;
f=p_2^{-1}e^{p/m}-p_1.$$
\eT
\begin{pf} \par Let $m=\deg P$ and
$P(f)=P_mf^m+\dots+P_0=Re^p$, where $P_j\in\R^n$. Clearly
$m\geq1$, since $P(f)$ is a non-trivial invertible function. As
the order of $f$ is finite, $p\in\P^n$. If $p$ is a constant it
follows that $f$ is a polynomial. Hence $\deg p\geq1$ and the
order of $f$ is positive.
\par Let us assume at first that $n=1$.
There is $r_0>0$ such that the functions $P_j$ and $R$ have
no zeros or poles when $|z|>r_0$. If $Q=R/P_m$ then $Qe^p$ is a
function of completely regular growth of order $\deg p$. By
Corollary \ref{C:crg}, $P(w)=P_m(w+p_1)^m$. Thus
$mP_mp_1=P_{m-1}$, so $p_1\in\R$. It follows that $P_m=Rp_2^m$,
$p_2\in\R$, and $f=p_2^{-1}e^{p/m}-p_1$.
\par If $n>1$ we let $z=(z',z_n)$, $z'\in{\Bbb C}^{n-1}$. There
exists a point $z_0\in{\Bbb C}^n$ so that the functions $P_m,\,R$
are holomorphic and non-vanishing in a neighborhood of $z_0$, and
the gradient of $p$ is also non-vanishing there. Without loss of
generality, we may assume that there is $\eps>0$ such that for any
fixed $z'$, $|z'|<\eps$, the polynomial $p(z',z_n)$ is not
constant and the functions $P_m(z',z_n)$ and $R(z',z_n)$ are well
defined rational functions not identically equal to 0.
\par From the one-dimensional case we see that
$$P(w)=P_m(z',z_n)(w+p_1(z',z_n))^m,$$ for each $z'$, $|z'|<\eps$.
Thus $mP_mp_1=P_{m-1}$, so $p_1\in\R^n$. Hence $P(w)=P_m(w+p_1)^m$
on an open set in ${\Bbb C}^n$ and, consequently, this equality
holds everywhere. We conclude that $P_m=Rp_2^m$ and
$f=p_2^{-1}e^{p/m}-p_1$, with $p_2\in\R^n$.
\end{pf}
\par This result allows us to describe all algebras $\R^n[f]$ with
non-trivial invertible elements in $\M_0^n$.
\bC\label{C:eq} Let
$f$ be an entire transcendental function on ${\Bbb C}^n$ of finite
order. The following statements are equivalent: \par (i) The
algebra $\R^n[f]$ contains non-trivial invertible elements in
$\M_0^n$. \par (ii) The function $f=q_1e^p+q_2$, where
$q_1,q_2\in\R^n$ and $p\in\P^n$.
\par (iii) $\R^n[f]=\R^n[e^p]$, where $p\in\P^n$. \par Moreover,
if one of these statements holds and $g\in\R^n[f]$ is invertible
in $\M_0^n$, then $g=re^{mp}$, where $m\ge0$ is an integer and
$r\in\R^n$. \eC
\begin{pf} We have $(i)$ implies $(ii)$ by Theorem \ref{T:fn0}. Clearly,
$(ii)$ implies $(iii)$, and $(iii)$ implies $(i)$.
\par If one of these statements holds and $g\in\R^n[f]$ is invertible
in $\M^n_0$, then $g=re^q$, where $r\in\R^n$ and $q\in\P^n$. By
Theorem \ref{T:fn0}, $q=mp$ for some integer $m\ge0$.
\end{pf}

We will need the following lemma, for the proof of our division
theorems.

\bL\label{L:Bezout} Let $f$ be an entire function on ${\Bbb C}^n$
and $P,Q$ be relatively prime polynomials of $\R^n[w]$. If the
function $P(f)/Q(f)\in\M_0^n$, then $Q(f)$ is an invertible
element in $\M_0^n$.\eL

\begin{pf} It suffices to show that the zero variety $Z(Q(f))$ of
the meromorphic function $Q(f)$ is algebraic. Assuming the
contrary, there exists an irreducible component $X$ of $Z(Q(f))$
which is not algebraic. Since $P(f)/Q(f)\in\M_0^n$, it follows
that $X\subset Z(P(f))$. As $P,Q$ are relatively prime, there
exist polynomials $P_1,Q_1\in\R^n[w]$ so that $PP_1+QQ_1=1$, hence
$P(f)P_1(f)+Q(f)Q_1(f)=1$. Note that the polar varieties of the
meromorphic functions $P_1(f),Q_1(f)$ are algebraic. As $X$ is not
algebraic, we can find a point $z_0\in X$ so that the functions
$P(f),P_1(f),Q(f),Q_1(f)$ are holomorphic in a neighborhood of
$z_0$. Since $P(f),Q(f)$ vanish at $z_0$, we get a
contradiction.\end{pf}
\par Let $\A$ be an algebra, $\B\subset\A$ a subalgebra and $g\in\B$.
It is an interesting problem to compare the ideal $\B g$ generated
by $g$ in $\B$ to the ideal $\A g\cap\B$. Clearly, $\B
g\subseteq\A g\cap\B$, but the inclusion is in general strict. If
$\A=\M_0^n$ and $\B=\R^n[f]$ we are going to describe
$\M_0^ng\cap\R^n[f]$ completely.
\par Given a function $g\in\R^n[e^p]$, where $p\in\P^n$ is not
constant, then $g=P(e^p)$ for a unique $P\in\R^n[w]$. Let $m\ge0$
denote the order of the zero of $P(w)$ at $w=0$. Then $g$ can be
written uniquely in the form $g=Q(e^p)e^{mp}$, where $m\ge0$,
$Q\in\R^n[w]$, $Q(0)\neq0$.

\bT\label{T:mdt} (i) Let $g\in\R^n[e^p]$, where $p\in\P^n$ is not
constant. If $g=Q(e^p)e^{mp}$, where $m\ge0$, $Q\in\R^n[w]$,
$Q(0)\neq0$, then
$$\M_0^ng\cap\R^n[e^p]=\R^n[e^p]Q(e^p).$$
\par (ii) Let $g\in\R^n[f]$, where $f$ is an entire function on
${\Bbb C}^n$ such that $\R^n[f]$ does not contain non-trivial
invertible elements in $\M_0^n$. Then
$$\M_0^ng\cap\R^n[f]=\R^n[f]g.$$
\eT
\begin{pf} $(i)$ Clearly, $Q(e^p)\in\M^n_0g$, so
$\R^n[e^p]Q(e^p)\sbs\M_0^ng\cap\R^n[e^p]$.
\par Suppose $h\in\M_0^n$ and $hg=P_0(e^p)$, where $P_0\in\R^n[w]$.
Let $P\in\R^n[w]$ be the greatest common divisor of $P_0,Q$, and
write $P_0=PQ_0$, $Q=PQ_1$, where $Q_j\in\R^n[w]$, $j=0,1$. Then
$he^{mp}=Q_0(e^p)/Q_1(e^p)\in\M_0^n$, so by Lemma \ref{L:Bezout},
$Q_1(e^p)$ is invertible in $\M^n_0$. Corollary \ref{C:eq} implies
that $Q_1(e^p)=re^{dp}$, where $r\in\R^n$ and $d\geq0$. Since
$Q(0)\neq0$ it follows that $d=0$. Hence $P_0=QQ_2$, for some
$Q_2\in\R^n[w]$, and $hg=Q_2(e^p)Q(e^p)$. So we have shown that
$\M_0^ng\cap\R^n[e^p]\sbs\R^n[e^p]Q(e^p)$.

\par $(ii)$ We have $\R^n[f]g\sbs\M^n_0g\cap\R^n[f]$. The opposite
inclusion is clearly true if $f$ is a polynomial, so we assume
that $f$ is transcendental.
\par Suppose that $g=P_1(f)$ and $hg=P_0(f)$ for some
$h\in\M_0^n$, where $P_0,P_1\in\R^n[w]$. We write as in $(i)$
$P_j=PQ_j$, where $P,Q_j\in\R^n[w]$, $j=0,1$, and $Q_0,Q_1$ are
relatively prime. Then we conclude by Lemma \ref{L:Bezout} that
$Q_1(f)$ is invertible in $\M^n_0$. Since $\R^n[f]$ does not
contain non-trivial invertible elements in $\M_0^n$ and since $f$
is transcendental, it follows that $\deg Q_1=0$. So $P_0=P_1Q_2$,
where $Q_2\in\R^n[w]$. Hence $hg=Q_2(f)g\in\R^n[f]g$.
\end{pf}
\par The immediate consequence of Theorem \ref{T:mdt} is
the division theorem for algebras $\R^n[f]$, which finishes the
proof of Theorem \ref{T:it}.

\bC\label{C:dt} Let $f$ be an entire function on ${\Bbb C}^n$,
$h_1,h_0\in\R^n[f]$, and assume that $g=h_0/h_1\in\M_0^n$.
\par (i) If $\R^n[f]$ does not contain non-trivial invertible
elements in $\M_0^n$, then $g\in\R^n[f]$.
\par (ii) If $\R^n[f]=\R^n[e^p]$, where $p\in\P^n$ is not constant,
then $g\in\R^n[e^p,e^{-p}]$. \eC

\par For functions $g_1,g_2\in\M^n_0$ we say that $g_1\sim g_2$ if
$g_1/g_2,g_2/g_1\in\M^n_0$. The following corollary  is an
immediate consequence of the previous results.
\bC\label{C:eqr} (i) Let $p\in\P^n$ be non-constant and
$g_1,g_2\in\R^n[e^p]$. Then $g_1\sim g_2$ if and only if
$g_2=Re^{mp}g_1$, for some $R\in\R^n$ and $m\in{\Bbb Z}$.
\par (ii) Let $f$ be an entire function on ${\Bbb C}^n$ and
$g_1,g_2\in\R^n[f]$. If $\R^n[f]$ does not contain non-trivial
invertible elements in $\M_0^n$, then $g_1\sim g_2$ if and only if
$g_2=Rg_1$, for some $R\in\R^n$. \eC

\section{Algebraic dependence over ${\Bbb C}$}\label{S:pairs}

\par Two meromorphic functions $f,\,g$ on ${\Bbb C}^n$ are
called algebraically dependent if there exists a nontrivial
polynomial $P$ on ${\Bbb C}^2$ so that $P(f,g)=0$. We now give a
complete characterization of such pairs of functions.

\begin{Theorem}\label{T:entpair} The entire functions
$f,g$ on ${\Bbb C}^n$ are algebraically dependent if and only if
there exists an entire function $h$ on ${\Bbb C}^n$ such that one
of the following holds:

\par (i) $f,g\in{\Bbb C}[h]$.

\par (ii) $f,g\in{\Bbb C}[e^h,e^{-h}]$.\end{Theorem}

\begin{pf} Two functions $f,g$ as in $(i)$ or $(ii)$ are clearly
algebraically dependent.  Conversely, let $P(z_1,z_2)$ be an
irreducible polynomial of degree $d$ so that $P(f,g)=0$. If $d=1$
then case $(i)$ in the statement clearly holds, so we may assume
that $d\geq2$. Moreover, we may assume that none of the maps $f,g$
is constant. Let $[z_0:z_1:z_2]$ denote the homogeneous
coordinates on the complex projective space ${\Bbb P}^2$, and
consider the standard embedding ${\Bbb C}^2\hookrightarrow{\Bbb
P}^2$, $(z_1,z_2)\to[1:z_1:z_2]$.

\par Let $X\subset{\Bbb P}^2$ be the algebraic curve defined by $P$,
$$X=\{[z_0:z_1:z_2]:\,\widetilde P(z_0,z_1,z_2)=0\},\;
\widetilde P(z_0,z_1,z_2)=z_0^dP(z_1/z_0,z_2/z_0).$$ We denote by
$X_\infty=X\cap\{z_0=0\}$ the set of points where $X$ intersects
the line at infinity. We then have a holomorphic mapping
$$F:{\Bbb C}^n\rightarrow X\setminus X_\infty,\;
F(\zeta)=[1:f(\zeta):g(\zeta)].$$ Consider the normalization of
$X$, $\sigma:\widetilde X\rightarrow X$, where $\widetilde X$ is a
compact Riemann surface (see e.g. \cite{Gr}), and let
$S=\sigma^{-1}(X_\infty)$. Finally, let
$\pi:Y\rightarrow\widetilde X$ be the universal covering of
$\widetilde X$, where $Y={\Bbb P}^1$, or $Y={\Bbb C}$, or
$Y=\Delta$ (the unit disk), and let $Z=\pi^{-1}(S)$. There exists
a non-constant holomorphic lifting of $F$, $G:{\Bbb
C}^n\rightarrow Y\setminus Z$, $F=\sigma\circ\pi\circ G$.
Therefore $Y\neq\Delta$. Suppose $Y={\Bbb C}$. Since
$S\neq\emptyset$ it follows that $Z$ is an infinite discrete
subset of ${\Bbb C}$, so $G$ is constant by Picard's theorem.

\par We conclude that $Y=\widetilde X={\Bbb P}^1$, and we have a
non-constant holomorphic map $G:{\Bbb C}^n\rightarrow {\Bbb P}^1
\setminus S$, $F=\sigma\circ G$. Picard's theorem implies that
$|S|\leq2$. Now $\sigma:{\Bbb P}^1\rightarrow{\Bbb P}^2$ is a
holomorphic map, with $\sigma({\Bbb P}^1)=X$. It follows that
$\sigma$ is a rational map, i.e. there exist homogeneous
polynomials $P_0,\,P_1,\,P_2$ of degree $m$ so that
$$\sigma([x_0:x_1])=[P_0(x_0,x_1):P_1(x_0,x_1):P_2(x_0,x_1)],
\;[x_0:x_1]\in{\Bbb P}^1.$$ We have two cases.

\par {\em Case 1.} $|S|=1$. Composing with a M\"obius map, we may
assume that $S=\{[0:1]\}$. It follows that
$G(\zeta)=[1:h(\zeta)]$, for some entire function $h$ on ${\Bbb
C}^n$. Since $S=\sigma^{-1}(\{z_0=0\})$, we see that
$P_0(x_0,x_1)=0$ if and only if $x_0=0$, hence
$P_0(x_0,x_1)=x_0^m$. Therefore
$$F(\zeta)=\sigma([1:h(\zeta)])=
[1:P_1(1,h(\zeta)):P_2(1,h(\zeta))],$$ and we conclude that
$f,g$ verify case $(i)$ from the statement.

\par {\em Case 2.} $|S|=2$. We may now assume that
$S=\{[0:1],[1:0]\}$, so $G(\zeta)=[1:e^{h(\zeta)}]$, where $h$ is
an entire function on ${\Bbb C}^n$. As in Case 1 we see that
$P_0(x_0,x_1)=0$ if and only if $x_0=0$ or $x_1=0$, hence
$P_0(x_0,x_1)=x_0^kx_1^{m-k}$, for some $1\leq k<m$. This yields
$$F(\zeta)=[1:e^{-Nh(\zeta)}P_1(1,e^{h(\zeta)}):
e^{-Nh(\zeta)}P_2(1,e^{h(\zeta)})],$$ where $N=m-k>0$, so the
functions $f,g$ verify $(ii)$.\end{pf}

\begin{Theorem}\label{T:meropair} The meromorphic functions $f,g$
on ${\Bbb C}^n$ are algebraically dependent if and only if one of
the following holds:

\par (i) There exists a meromorphic function $h$ on ${\Bbb C}^n$
and rational functions $R_1,R_2$ on ${\Bbb C}$, so that
$f=R_1\circ h$, $g=R_2\circ h$.

\par (ii) There exists an entire function $h$ on ${\Bbb C}^n$,
and elliptic functions $\varphi_1,\varphi_2$ with the same
periods, so that $f=\varphi_1\circ h$, $g=\varphi_2\circ
h$.\end{Theorem}

\begin{pf} Two functions $f,g$ as in $(i)$ are clearly
algebraically dependent. If $f,g$ are as in $(ii)$, then
$\sigma=[1:\varphi_1:\varphi_2]$ induces a holomorphic map from a
complex torus into ${\Bbb P}^2$. The image of this map is an
analytic variety, hence algebraic by Chow's theorem. This implies
that $f,g$ are algebraically dependent.

\par Conversely, assume that $f,g$ are algebraically dependent and
non-constant, and let $P(z_1,z_2)$ be an irreducible polynomial of
degree $d$ so that $P(f,g)=0$. We use the same notation as in the
proof of Theorem \ref{T:entpair}: $X=\{P=0\}\subset{\Bbb P}^2$,
$\sigma:\widetilde X\rightarrow X$ is the normalization of $X$,
and $\pi:Y\rightarrow\widetilde X$ is the universal covering of
$\widetilde X$.

\par Let $A$ be the union of the polar varieties of $f$ and $g$.
We define the meromorphic map
$F=[f_1g_1:f_0g_1:g_0f_1]:D\dashrightarrow X$, where $f=f_0/f_1$,
$g=g_0/g_1$, and $f_j,\,g_j$ are entire functions. We have
$F=[1:f:g]$ on ${\Bbb C}^n\setminus A$. If $I$ is the
indeterminacy set of $F$, then $I$ is an analytic subvariety of
${\Bbb C}^n$ of codimension at least 2 and $I\subset A$. Hence
$D={\Bbb C}^n\setminus I$ is simply connected, so $F$ has a
non-constant holomorphic lifting $G:D\rightarrow Y$,
$F=\sigma\circ\pi\circ G$. If $Y=\Delta$, then, since ${\rm
codim}\,I\geq2$, $G$ extends holomorphically to ${\Bbb C}^n$, so
it is constant. Therefore we have two cases.

\par {\em Case 1.} $Y=\widetilde X={\Bbb P}^1$. Then
$F=\sigma\circ G$, $G=[h_0:h_1]:D\rightarrow {\Bbb P}^1$, where
$h_j$ are holomorphic on $D$, so they extend holomorphically to
${\Bbb C}^n$. Moreover $\sigma=[P_0:P_1:P_2]$, for some
homogeneous polynomials $P_0,\,P_1,\,P_2$ of degree $m$. It
follows that $f,g$ verify case $(i)$ from the statement, with
$R_1(x)=P_1(1,x)/P_0(1,x)$, $R_2(x)=P_2(1,x)/P_0(1,x)$,
$h=h_1/h_0$.

\par {\em Case 2.} $Y={\Bbb C}$, so $\widetilde X$ is a complex
torus. Then $G$ extends to an entire function on ${\Bbb C}^n$,
since ${\rm codim}\,I\geq2$. The holomorphic map
$\sigma:\widetilde X\rightarrow X\subset{\Bbb P}^2$ is given by
$\sigma=[\psi_0:\psi_1:\psi_2]$, where $\psi_j$ are meromorphic
functions on $\widetilde X$. It follows that
$\rho_j=\psi_j\circ\pi$ are elliptic functions with the same
periods, and $f,\,g$ verify $(ii)$ with $\varphi_1=\rho_1/\rho_0$,
$\varphi_2=\rho_2/\rho_0$.\end{pf}

\end{document}